\documentclass[12pt]{amsart}
\usepackage{amssymb,latexsym,comment,url}

\newcommand{\Fq}{{\mathbb{F}_q}}
\newcommand{\C}{{\mathbb{C}}}
\newcommand{\Fp}{{\mathbb{F}_p}}

\newcommand{\isom}{ \cong }

\newcommand{\Z}{{\mathbb Z}}

\newenvironment{Proof}{\par\noindent{\sc Proof:}}%
                      {\hspace*{\fill}\nobreak$\Box$\par\medskip}
                       {\hspace*{\fill}\nobreak$\Box$\par\medskip}

\newtheorem{Proposition}{Proposition}[section]
\newtheorem{Theorem}[Proposition]{Theorem}
\newtheorem{Lemma}[Proposition]{Lemma}
\newtheorem{Corollary}[Proposition]{Corollary}

\theoremstyle{definition}

\renewcommand{\baselinestretch}{1.1}

\begin{document}

\title[Hypergeometric Series and Huff Curves]%
{Evaluation of Gaussian hypergeometric series using Huff's models of elliptic curves}
\author[M. Sadek]%
{Mohammad~Sadek}
\address{American University in Cairo, Department of Mathematics and Actuarial Science, AUC Avenue, New Cairo, Egypt}
\email{mmsadek@aucegypt.edu}

\author[N. El-Sissi]%
{Nermine El-Sissi}
\address{American University in Cairo, Department of Mathematics and Actuarial Science, AUC Avenue, New Cairo, Egypt}
\email{nelsissi@aucegypt.edu}

\author[A. S. Zargar]%
{Arman Shamsi Zargar}
\address{Department of Mathematics and Applications, Faculty of Science, University of Mohaghegh Ardabili, Ardabil 56199--11367, Iran}
\email{shzargar.arman@gmail.com}

\author[N. Zamani]%
{Naser Zamani}
\address{Department of Mathematics and Applications, Faculty of Science, University of Mohaghegh Ardabili, Ardabil 56199--11367, Iran}
\email{naserzaka@yahoo.com}

\let\thefootnote\relax

\begin{abstract}{\footnotesize A Huff curve over a field $K$ is an elliptic curve defined by the equation $ax(y^2-1)=by(x^2-1)$ where $a,b\in K$ are such that $a^2\ne b^2$. In a similar fashion, a general Huff curve over $K$ is described by the equation $x(ay^2-1)=y(bx^2-1)$ where $a,b\in K$ are such that $ab(a-b)\ne 0$. In this note we express the number of rational points on these curves over a finite field $\Fq$ of odd characteristic in terms of Gaussian hypergeometric series $\displaystyle {_2F_1}(\lambda):={_2F_1}\left(\begin{matrix}   \phi&\phi \\  & \epsilon  \end{matrix}\Big| \lambda \right)$ where $\phi$ and $\epsilon$ are the quadratic and trivial characters over $\Fq$, respectively. Consequently, we exhibit the number of rational points on the elliptic curves $y^2=x(x+a)(x+b)$ over $\Fq$ in terms of ${_2F_1}(\lambda)$. This generalizes earlier known formulas for Legendre, Clausen and Edwards curves. Furthermore, using these expressions we display several transformations of ${_2F_1}$. Finally, we present the exact value of $_2F_1(\lambda)$ for different $\lambda$'s over a prime field $\Fp$ extending previous results of Greene and Ono.  }
\end{abstract}

\maketitle

\section{Introduction}
An elliptic curve $E$ over a field $K$ is a smooth projective curve of genus 1 equipped with a $K$-rational point. When embedded in the projective plane, $E$ can be described by a Weierstrass equation of the form
\[E: \ y^2+a_1xy+a_3y=x^3+a_2x^2+a_4x+a_6\]
where $a_1, a_2, a_3, a_4, a_6\in K$. Elliptic curves can be represented by several other equations. The interested reader may consult \cite[Chapter 2]{Wa}. In the few past decades, many alternative equations describing $E$ have been introduced due to cryptographic applications. In the cryptographic settings, curves are usually considered over finite fields $\Fq$ with $q$ elements.

Studying rational points of elliptic curves over $\Fq$ traces back to Gauss \cite[Page 68]{S.Z}, \cite[Page 2]{Sil}. In addition, finding the number of rational points of an elliptic curve over a given finite field $\Fq$ has always been of interest to many researchers. Schoof \cite{Schoof} gives three algorithms that compute the number of these rational points. This motivates cryptographers to present practical algorithms to enumerate points on higher genus curves, see \cite{Gaudry}.

Greene \cite{Greene} initiated the study of Gaussian hypergeometric series over finite fields. This was followed by the work of Ono attempting to link the number of rational points on elliptic curves described by Legendre equations over finite fields with certain values of these Gaussian hypergeometric series, see \cite{Ono1}. Influenced by the latter articles, research has been targeting the investigation of the various links between elliptic curves described by different equations and Gaussian hypergeometric series. For example, similar links were found for Clausen curves \cite{Guindy} and Edwards curves \cite{Sadek1}. Further connections between these series and the number of rational points on higher genus curves were explored, see \cite{BarmanKalita,Barmanalgebraiccurves,BarmanKalitaSaikia,Sadek2}

Herein, we consider Huff curves which are elliptic curves defined by an equation of the form
\[H_{a,b}: ax(y^2-1)=by(x^2-1), \;a,b\in\mathbb{F}_q^{\times},\;a^2\ne b^2.\]
 These curves were first studied by Huff to answer a Diophantine question on rational distance sets, see \cite{Huff}. Huff curves were recently revisited due to the completeness of the addition law on them. This means that the formulas for the addition work for all pairs of points on Huff curves with no exceptions. This turns out to be useful for cryptographic applications, see \cite{Joye}. In \cite{Wu}, general Huff curves were introduced. They are elliptic curves defined by the equation
\[G_{a,b}: x(ay^2-1)=y(bx^2-1),\; ab(a-b)\ne 0.\]
We show that the number of rational points on Huff and general Huff curves over $\Fq$ can be written in terms of a Gaussian hypergeometric series of type ${_2F_1}$. Making use of explicit isomorphisms between Huff curves and elliptic curves with complete $2$-torsion, we are able to describe the number of $\Fq$-rational points on elliptic curves described by the Weierstrass equation $y^2=x(x+a)(x+b)$ using Gaussian hypergeometric series $_2F_1$. This generalizes similar results of Ono on Legendre elliptic curves \cite{Ono1}.

Exploiting the existence of an isomorphism between Huff curves and Edwards curves, we are able to find some transformations of Gaussian hypergeometric series $\displaystyle {_2F_1}(\lambda):={_2F_1}\left(\begin{matrix}   \phi&\phi \\  & \epsilon  \end{matrix}\Big| \lambda \right)$ using the counts in \cite{Sadek1}. Additionally, these transformations are deployed to evaluate $_2F_1(\lambda)$ at new values of $\lambda$ extending on the findings of \cite{Greene} and \cite{Ono1}.

\section{Preliminaries}

In this note, $p$ will be an odd prime and $q=p^r$ for some positive integer $r$. A multiplicative character $\chi$ over $\mathbb{F}_q^{\times}$ will be extended to be defined over $\mathbb{F}_q$ by setting $\chi(0)=0$. We will write $\overline{\chi}$ for $1/\chi$. The Greek letters $\epsilon$ and $\phi$ are preserved for the trivial and quadratic characters respectively. More precisely, $\epsilon(x)=1$ for every $x\in\mathbb{F}_q^{\times}$, and $\epsilon(0)=0$; whereas $\phi(x)=1$ if $x$ is a square in $\mathbb{F}_q^{\times}$, $\phi(x)=-1$ if $x$ is not a square, and $\phi(0)=0$.

Let $J(A,B)$ be the Jacobi sum \[J(A,B)=\sum_{x\in\Fq}A(x)B(1-x)\] where $A$ and $B$ are multiplicative characters over $\Fq$. Let $A_0,A_1,\ldots,A_n$ and $B_1,\ldots,B_n$ be multiplicative characters defined over $\mathbb{F}_q$. The {\em Gaussian hypergeometric series} $_{n+1}F_n$ is defined as follows \[_{n+1}F_n\left(\begin{matrix}
                          A_0&A_1&\ldots&A_n \\
                          {} & B_1&\ldots&B_n
                        \end{matrix}\Big| x\right):=\frac{q}{q-1}\sum_{\chi}{{A_0\chi}\choose{\chi}}{{A_1\chi}\choose{B_1\chi}}\ldots{{A_n\chi}\choose{B_n\chi}}\chi(x)\]
where the sum is over all multiplicative characters over $\mathbb{F}_q$ and \[{A\choose B}:=\frac{B(-1)}{q}J(A,\overline{B})=\frac{B(-1)}{q}\sum_{x\in\mathbb{F}_q}A(x)\overline{B}(1-x).\]

The following lemma contains some of the properties of the symbol $\displaystyle{A\choose B}$ that can be found in \cite{Greene}.

\begin{Lemma}
\label{lem:charactersproperties}
Let $\epsilon$ and $\phi$ be the trivial and quadratic characters over $\mathbb{F}_q$ respectively. For any multiplicative characters $A$ and $B$ over $\mathbb{F}_q$, one has:
\begin{itemize}
\item[a)] $\displaystyle A(1+x)=\delta(x)+\frac{q}{q-1}\sum_{\chi}{A\choose \chi}\chi(x)$ where $\delta(x)=1$ if $x=0$ and $\delta(x)=0$ if $x\ne 0$,
\item[b)] $\displaystyle \overline{A}(1-x)=\delta(x)+\frac{q}{q-1}\sum_{\chi}{A\chi\choose \chi}\chi(x)$ where $\delta(x)=1$ if $x=0$ and $\delta(x)=0$ if $x\ne 0$,
\item[c)] $\displaystyle {A\choose B}={A\choose A\overline{B}}$,
\item[d)] $\displaystyle {A\choose B}={B\overline{A}\choose B}B(-1)$,
\item[e)] $\displaystyle {A\choose \epsilon}={A\choose A}=-\frac{1}{q}+\frac{q-1}{q}\delta(A)$ where $\delta(A)=1$ if $A=\epsilon$ and $\delta(A)=0$ otherwise,
\item[f)] $\displaystyle {B^2\chi^2\choose \chi}={\phi B \chi\choose\chi}{B\chi\choose B^2\chi}{\phi\choose\phi B}^{-1}B\chi(4)$,
\item[g)] $\displaystyle {A^2\choose AB}={A\choose B}{\phi A\choose A B} {\phi \choose B}^{-1}A(4)$.
\end{itemize}
\end{Lemma}

\section{Evaluation of some character sums}
In this section, we evaluate the character sums that we use in this work.

\begin{Lemma}
\label{lem1}
Let $f:\Fq\to\C$ be a map. One has
\[\sum_{x\in\mathbb{F}_q}\phi(x)f(x)=\sum_{x\in\mathbb{F}_q}f(x^2)-\sum_{x\in\mathbb{F}_q}f(x).\]
\end{Lemma}
\begin{Proof}
This holds since
\begin{eqnarray*}
\sum_{x\in\mathbb{F}_q}\phi(x)f(x)=\sum_{x\in\mathbb{F}_q}(1+\phi(x))f(x)-\sum_{x\in\mathbb{F}_q}f(x).
\end{eqnarray*}
This follows because $1+\phi(x)$ is either $2$ if $x$ is a square in $\mathbb{F}_q^{\times}$; $1$ if $x=0$; or $0$ otherwise.
\end{Proof}
The following corollary can be found in \cite{Sadek2}. We reproduce the proof for the convenience of the reader.
\begin{Corollary}
\label{cor:charactersum}
Let $a\in\mathbb{F}^{\times}_q$ and let $\chi, \psi$ be characters on $\mathbb{F}_q$. Then
\[\sum_{x\in\mathbb{F}_q}\psi(x^2)\chi(1+ax^2)=\psi\left(-a^{-1}\right)J(\psi,\chi)+\phi\psi\left(-a^{-1}\right)J(\phi\psi,\chi).\]
In particular, if $\chi=\phi$ then
\[\sum_{x\in\mathbb{F}_q}\psi(x^2)\phi(1+ax^2)=q\phi(-1)\left[\psi\left(-a^{-1}\right){\psi\choose\phi\psi}+\phi\psi\left(-a^{-1}\right){\phi\psi\choose\psi}\right].\]
\end{Corollary}
\begin{Proof}
By setting $f(x)=\psi(x)\chi(1+ax)$, Lemma \ref{lem1} implies that
$$\sum_{x\in\Fq}\psi(x^2)\chi(1+ax^2)
=\sum_{x\in\Fq}\psi(x)\chi(1+ax)+\sum_{x\in\Fq}\phi\psi(x)\chi(1+ax).$$
Now, since the map $x\mapsto -a^{-1}x$ is bijective over $\mathbb{F}_q$, one has
\begin{align*}
\sum_{x\in\Fq}\psi(x^2)\chi(1+ax^2)&=\psi\left(-a^{-1}\right)\sum_{x\in\Fq}\psi(x)\chi(1-x)+\phi\psi\left(-a^{-1}\right)\sum_{x\in\Fq}\phi\psi(x)\chi(1-x)\\
&=\psi\left(-a^{-1}\right)J(\psi,\chi)+\phi\psi\left(-a^{-1}\right)J(\phi\psi,\chi)
\end{align*}
where the last equality holds using the definition of the Jacobi sum.
Putting $\chi=\phi$, one has
\begin{align*}
\sum_x\psi(x^2)\phi(1+ax^2)&=\psi\left(-a^{-1}\right)J(\psi,\phi)+\phi\psi\left(-a^{-1}\right)J(\phi\psi,\phi)\\
&=q\phi(-1)\left[\psi\left(-a^{-1}\right){\psi\choose\phi}+\phi\psi\left(-a^{-1}\right){\phi\psi\choose\phi}\right]
\end{align*}
where the last equality is implied by the definition of the symbol $\displaystyle{A\choose B}$.
Now one concludes using Lemma \ref{lem:charactersproperties} (c).
\end{Proof}

\section{Number of rational points on Huff curves}
The following lemma evaluates certain symbols $\displaystyle  {A\choose B}$.
\begin{Lemma}
\label{lem2}
For any multiplicative character $\chi$ over $\Fq$, one has
\begin{itemize}
\item[a)] $\displaystyle {\phi\choose\chi}={{\phi\chi}\choose \chi} \chi(-1)$,
\item[b)] $\displaystyle {\epsilon \choose\epsilon}=\frac{q-2}{q}$,
\item[c)] $\displaystyle{\chi^2\choose\chi}={{\phi\chi}\choose\chi} \chi(4)$ if $\chi\ne\epsilon$.
\end{itemize}
\end{Lemma}
\begin{Proof}
The first equality follows directly from Lemma \ref{lem:charactersproperties} (d) while the second equality follows from Lemma \ref{lem:charactersproperties} (e). By setting $B=\epsilon$ in Lemma \ref{lem:charactersproperties} (f), we get
\begin{eqnarray*}
{\chi^2\choose\chi}={{\phi\chi}\choose\chi}{\chi\choose\chi}{\phi\choose\phi}^{-1}\chi(4).
\end{eqnarray*}
Now Lemma 2.1 (e) implies that $\displaystyle {\chi\choose\chi}={\phi\choose\phi}=-1/q$, if $\chi\ne\epsilon$.
\end{Proof}

A {\em general Huff curve} is defined by the following equation
\[G_{a,b}:x(ay^2-1)=y(bx^2-1),\; ab(a-b)\ne 0.\]
 We notice first that the curve has three points at infinity. This follows by homogenizing the defining equation via $x=X/Z$, $y=Y/Z$, so that we get the general Huff curve in the projective plane as
$X(aY^2-Z^2)=Y(bX^2-Z^2)$. Hence, the points $(1:0:0)$, $(0:1:0)$ and $(a:b:0)$ are in $G_{a,b}(\Fq)$.

The following theorem gives the number of rational points on $G_{a,b}$ over $\Fq$ in terms of a hypergeometric series of type ${_2F_1}$.
\begin{Theorem}
\label{Thm:rationalpointsongeneralHuffcurves}
Let $G_{a,b}$ be a general Huff curve defined by the equation
$x(ay^2-1)=y(bx^2-1)$ where $ab(a-b)\ne 0$. Then
\[|G_{a,b}(\Fq)|=q+2-\frac{1}{q-1}-\left(2+\frac{1}{q-1}\right)\phi(a^{-1}b)+\frac{q^2}{q-1}{_2F_1}\left(\begin{matrix}
                          \phi&\phi \\
                          {} & \epsilon
                        \end{matrix}\Big| a^{-1}b\right).\]
\end{Theorem}
\begin{Proof}
After completing the square, the curve $G_{a,b}$ is described by the equation
\[Y^2=\frac{b^2x^4+(4a-2b)x^2+1}{4a^2x^2},\quad x\ne 0,\]
where $\displaystyle Y=y-\frac{bx^2-1}{2ax}$, $x\ne 0$.
In fact, there is only one point in $G_{a,b}(\Fq)$ that has a zero $x$-coordinate. Therefore, we obtain the following equalities
\begin{align*}
|G_{a,b}(\Fq)|&=q+3+1+\sum_{x\in\mathbb{F}_q^{\times}}\phi\left(b^2x^4+(4a-2b)x^2+1\right)\\
&=q+4+\sum_{x\in\mathbb{F}_q^{\times}}\phi\left(4ax^2+(1-bx^2)^2\right)\\
&=q+4+\sum_{x\in\mathbb{F}_q^{\times}}\phi\left(1+\frac{a^{-1}}{4x^2}(1-bx^2)^2\right).
\end{align*}
Lemma 2.1 (a) gives the following equalities
\begin{align*}
|G_{a,b}(\Fq)|&=q+4+\sum_{x\in\mathbb{F}_q^{\times}}\left\{\frac{q}{q-1}\sum_{\chi}{\phi\choose\chi}\chi\left(\frac{a^{-1}}{4x^2}(1-bx^2)^2\right)\right\}\\
&=q+4+\frac{q}{q-1}\sum_{x\in\mathbb{F}_q^{\times}}\sum_{\chi}{\phi\choose\chi}\chi\left(\frac{a^{-1}}{4x^2}\right)\chi^2(1-bx^2)\\
&=q+4+\frac{q}{q-1}\sum_{\chi}{\phi\choose\chi}\chi\left(\frac{a^{-1}}{4}\right)\left(\sum_{x\in\mathbb{F}_q^{\times}}\bar{\chi}(x^2)\chi^2(1-bx^2)\right).
\end{align*}
On account of Corollary \ref{cor:charactersum}, we get
\begin{align*}
|G_{a,b}(\Fq)|&=q+4+\frac{q}{q-1}\sum_{\chi}{\phi\choose\chi}\chi\left(\frac{a^{-1}}{4}\right)\left[\bar\chi(b^{-1})J(\bar{\chi},\chi^2)+\phi\bar\chi(b^{-1})J(\phi\bar{\chi},\chi^2)\right]\\
&=q+4+\frac{q}{q-1}\sum_{\chi}{\phi\choose\chi}\chi\left(\frac{a^{-1}}{4}\right)\chi(b)J(\bar\chi,\chi^2)+\frac{q}{q-1}\sum_{\chi}{\phi\choose\chi}\chi\left(\frac{a^{-1}}{4}\right)\phi\chi(b)J(\phi\bar\chi,\chi^2).
\end{align*}
Using the identities
$$J(\bar\chi,\chi^2)=J(\chi^2,\bar\chi)=q\bar\chi(-1){\chi^2\choose\chi}, \ J(\phi\bar\chi,\chi^2)=J(\chi^2,\phi\bar\chi)=q\phi\bar\chi(-1){\chi^2\choose\phi\chi},$$
we get
\begin{align*}
|G_{a,b}(\Fq)|
&=q+4+\frac{q^2}{q-1}\sum_{\chi}{\phi\choose\chi}{\chi^2\choose\chi}\chi\left(-\frac{a^{-1}b}{4}\right) +\frac{q^2}{q-1}\phi(-1)\sum_{\chi}{\phi\choose\chi}{\chi^2\choose\phi\chi}\chi\left(-\frac{a^{-1}b}{4}\right).
\end{align*}
We note that over $\mathbb{F}_q^{\times}$, we have $\bar\phi=\phi$ and $\phi^2=1$.

Now Lemma \ref{lem2} implies that
\begin{align*}
|G_{a,b}(\Fq)|&=q+4+\frac{q^2}{q-1}\sum_{\chi}{{\phi\chi}\choose\chi}{\chi^2\choose\chi}\chi\left(\frac{a^{-1}b}{4}\right) +\frac{q^2}{q-1}\phi(-1)\sum_{\chi}{{\phi\chi}\choose\chi}{\chi^2\choose\phi\chi}\chi\left(\frac{a^{-1}b}{4}\right)\\
&=q+4+\frac{q^2}{q-1}\sum_{\chi\ne\epsilon}{{\phi\chi}\choose\chi}{\chi^2\choose\chi}\chi\left(\frac{a^{-1}b}{4}\right)+\frac{q^2}{q-1}{\phi\choose\epsilon}{\epsilon\choose\epsilon} \\
&\qquad +\frac{q^2}{q-1}\phi(-1)\sum_{\chi}{{\phi\chi}\choose\chi}{\chi^2\choose\phi\chi}\chi\left(\frac{a^{-1}b}{4}\right).
\end{align*}
According to Lemma \ref{lem2} and Lemma \ref{lem:charactersproperties} (g), (e), the latter equality can be written as follows
\begin{align*}
  |G_{a,b}(\Fq)|&=q+4+\frac{q^2}{q-1}\sum_{\chi\ne\epsilon}{{\phi\chi}\choose\chi}{\phi\chi\choose\chi}\chi\left(a^{-1}b\right)-\frac{q-2}{q-1}\\
  &\qquad +\frac{q^2}{q-1}\phi(-1)\sum_{\chi}{{\phi\chi}\choose\chi}{\chi\choose\phi}{{\phi\chi}\choose{\phi\chi}}{\phi\choose\phi}^{-1}\chi\left(a^{-1}b\right)\\
  &=q+3+\frac{1}{q-1}+\frac{q^2}{q-1}\sum_{\chi\ne\epsilon}{{\phi\chi}\choose\chi}{\phi\chi\choose\chi}\chi\left(a^{-1}b\right)
   -\frac{q^3}{q-1}\phi(-1)\sum_{\chi}{{\phi\chi}\choose\chi}{\chi\choose\phi}{{\phi\chi}\choose{\phi\chi}}\chi\left(a^{-1}b\right)\\
   &=q+3+\frac{1}{q-1}+\frac{q^2}{q-1}\sum_{\chi}{{\phi\chi}\choose\chi}{\phi\chi\choose\chi}\chi\left(a^{-1}b\right)-\frac{q^2}{q-1}{\phi\choose\epsilon}^2\\
   &\qquad -\frac{q^3}{q-1}\phi(-1)\sum_{\chi\ne\phi}{{\phi\chi}\choose\chi}{\chi\choose\phi}{{\phi\chi}\choose{\phi\chi}}\chi\left(a^{-1}b\right)-\frac{q^3}{q-1}\phi(-1){\epsilon\choose\phi}{\phi\choose\phi}{\epsilon\choose\epsilon}\phi(a^{-1}b)\\
   &=q+3+\frac{q^2}{q-1}\sum_{\chi}{{\phi\chi}\choose\chi}{\phi\chi\choose\chi}\chi\left(a^{-1}b\right) +\frac{q^2}{q-1}\phi(-1)\sum_{\chi\ne\phi}{{\phi\chi}\choose\chi}{\chi\choose\phi}\chi\left(a^{-1}b\right)-\frac{q-2}{q-1}\phi(a^{-1}b)\\
   &=q+3-\frac{q-2}{q-1}\phi(a^{-1}b)+\frac{q^2}{q-1}\sum_{\chi}{{\phi\chi}\choose\chi}{\phi\chi\choose\chi}\chi\left(a^{-1}b\right) \\ &\qquad +\frac{q^2}{q-1}\phi(-1)\sum_{\chi}{{\phi\chi}\choose\chi}{\chi\choose\phi}\chi\left(a^{-1}b\right)-\frac{q^2}{q-1}\phi(-1){\epsilon\choose\phi}{\phi\choose\phi}\phi(a^{-1}b)\\
   &= q+3-\phi(a^{-1}b)+\frac{q^2}{q-1}\sum_{\chi}{{\phi\chi}\choose\chi}{\phi\chi\choose\chi}\chi\left(a^{-1}b\right)   +\frac{q^2}{q-1}\phi(-1)\sum_{\chi}{{\phi\chi}\choose\chi}{\chi\choose\phi\chi}\chi\left(a^{-1}b\right)\\
   &= q+3-\phi(a^{-1}b)+\frac{q^2}{q-1}\left[{_2F_1}\left(\begin{matrix}
                          \phi&\phi \\
                          {} & \epsilon
                        \end{matrix}\Big| a^{-1}b\right)+\phi(-1)\,{_2F_1}\left(\begin{matrix}
                          \phi&\epsilon \\
                          {} & \phi
                        \end{matrix}\Big| a^{-1}b\right)\right]
\end{align*}
where
\[{_2F_1}\left(\begin{matrix}
                          \phi&\epsilon \\
                          {} & \phi
                        \end{matrix}\Big| a^{-1}b\right)=-\frac{1}{q}\phi(-1)\left(1+\phi(a^{-1}b)\right),\]
                        see Corollary 3.16 of \cite{Greene}. Therefore, one may deduce that
                        \[|G_{a,b}(\Fq)|=q+3-\phi(a^{-1}b)-\frac{q}{q-1}\left(1+\phi(a^{-1}b)\right)+\frac{q^2}{q-1}{_2F_1}\left(\begin{matrix}
                          \phi&\phi \\
                          {} & \epsilon
                        \end{matrix}\Big| a^{-1}b\right).\]
\end{Proof}

Let $H_{a,b}$ be the curve defined by the following equation over $\Fq$
\[H_{a,b}: ax(y^2-1)=by(x^2-1), \;a,b\in\mathbb{F}_q^{\times},\;a^2\ne b^2.\]
We recall that the general Huff curve $G_{a,b}$ becomes a Huff curve when $a$ and $b$ are both squares in $\Fq$, see \cite[\S 1]{Wu}. More precisely, the curves $G_{a^2,b^2}$ and $H_{a,b}$ are isomorphic over $\Fq$. In particular, the family of Huff curves is contained in general Huff curves. Thus, we obtain the following result.

\begin{Corollary}
\label{cor:rationalpointsonHuffcurves}
Let $H_{a,b}$ be a Huff curve defined by the equation
$ax(y^2-1)=by(x^2-1)$ where $a,b\in\mathbb{F}_q^{\times}$ are such that $a^2\ne b^2$. Then
\[|H_{a,b}(\Fq)|=q-\frac{2}{q-1}+\frac{q^2}{q-1}\,{_2F_1}\left(\begin{matrix}
                          \phi&\phi \\
                          {} & \epsilon
                        \end{matrix}\Big| a^{-2}b^2\right).\]
\end{Corollary}
\begin{Proof}
This follows directly from Theorem \ref{Thm:rationalpointsongeneralHuffcurves} together with the observation that $H_{a,b}$ and $G_{a^2,b^2}$ are isomorphic over $\Fq$.
\end{Proof}

\section{Transformations and values of Gaussian hypergeometric series}
In this section we introduce several applications of the expressions in Theorem \ref{Thm:rationalpointsongeneralHuffcurves} and Corollary \ref{cor:rationalpointsonHuffcurves}. We firstly show that we managed to reach an expression for the number of rational points on elliptic curves defined by $y^2=x(x+a)(x+b)$ over $\Fq$ in terms of the Gaussian hypergeometric series $_2F_1$. Writing $_2F_1(\lambda)$ for $\displaystyle{_2F_1}\left(\begin{matrix}
                          \phi&\phi \\
                          {} & \epsilon
                        \end{matrix}\Big| \lambda\right)$, we then find a transformation between $_2F_1(\lambda^2)$, $\displaystyle _2F_1\left(\left(\frac{1-\lambda}{1+\lambda}\right)^2\right)$, and $\displaystyle _2F_1\left(\frac{4\lambda}{(1+\lambda)^2}\right)$. Finally, we show how we may find new values for $_2F_1(\lambda)$.

 The Huff curve $H_{a,b}: ax(y^2-1)=by(x^2-1)$ is isomorphic to the elliptic curve defined by $v^2=u(u+a^2)(u+b^2)$. Moreover, $E$ contains a subgroup isomorphic to $\Z/2\Z\times\Z/4\Z$. Conversely, any elliptic curve with rational subgroup isomorphic to $\Z/2\Z\times\Z/4\Z$ is isomorphic to a Huff curve. More specifically, $E$ contains a rational subgroup isomorphic to $\Z/2\Z\times\Z/4\Z$ if and only if it admits a Weierstrass equation of the form $y^2=x(x+t^2)(x+(t+2)^2)$, see \cite[Theorem 2]{Joye}.
In addition, a general Huff curves $G_{a,b}$ is isomorphic to an elliptic curve defined by the Weierstrass equation $v^2=u(u+a)(u+b)$ via the following transformations 
\begin{eqnarray*}
u= \frac{bx - ay}{y-x},\quad v=\frac{b-a}{y-x};\textrm{ whereas }x=\frac{u+a}{v},\quad y=\frac{u+b}{v}.
\end{eqnarray*}
For a full description of the isomorphisms above, the reader may consult \cite{Wu}.

 \begin{Corollary}
 Let $E_{a,b}$ be an elliptic curve defined by $y^2=x(x+a)(x+b)$ over $\Fq$. Then the number of rational points in $E(\Fq)$ is given as follows
 \[|E_{a,b}(\Fq)|=q+2-\frac{1}{q-1}-\left(2+\frac{1}{q-1}\right)\phi(a^{-1}b)+\frac{q^2}{q-1}\;{_2F_1}\left(\begin{matrix}
                          \phi&\phi \\
                          {} & \epsilon
                        \end{matrix}\Big| a^{-1}b\right).\]
 \end{Corollary}
 \begin{Proof}
 This is given by making use of the isomorphism between the Huff curve $G_{a,b}$ and $E_{a,b}$, see Theorem 2.1 in \cite{Wu} for an explicit description of the isomorphism. Now we use Theorem \ref{Thm:rationalpointsongeneralHuffcurves}.
 \end{Proof}
 Furthermore, the Huff curve $H_{a,b}$ is isomorphic to an Edwards curve described by $E_{d^2}:x^2+y^2=1+d^2x^2y^2$ where $\displaystyle d=\frac{a-b}{a+b}$, see \cite[\S 1]{Wu}. In \cite{Sadek1}, the size of the rational points on an Edwards curves is expressed in terms of $_2F_1$ over $\Fp$. Yet, all the results there are valid over $\Fq$. Thus, one gets the following.
 \begin{Theorem}
 \label{thm:transformation}
 Let $\lambda\ne0,\pm1$ be in $\Fq$. The following identities hold
 \begin{align*}
 {_2F_1}(\lambda^2)&=\frac{q+1}{q^2}+\frac{q-1}{q}\phi(-1)\;{_2F_1}\left(\left(\frac{1-\lambda}{1+\lambda}\right)^2\right)\\
 &=\frac{q+1}{q^2}+\frac{q-1}{q}\;{_2F_1}\left(\frac{4\lambda}{(1+\lambda)^2}\right)\\
 &=\frac{q+1}{q^2}+\frac{q-1}{q}\phi(\lambda)\;{_2F_1}\left(\frac{(1-\lambda)^2}{-4\lambda}\right).
 \end{align*}
 \end{Theorem}
\begin{Proof}
  One knows that $H_{a,b}$ is isomorphic to $E_{d^2}$ over $\Fq$ where $d=(a-b)/(a+b)$. In particular, $H_{1,\lambda}\isom E_{d^2}$ where $\displaystyle d=\frac{1-\lambda}{1+\lambda}$. One obtains the following equality by using Corollary \ref{cor:rationalpointsonHuffcurves} for $|H_{1,\lambda}(\Fq)|$ and \cite[Theorem 4.1]{Sadek1} for $|E_{d^2}(\Fq)|$
  \begin{align*}
  q-\frac{2}{q-1}+\frac{q^2}{q-1}\,{_2F_1}\left(\begin{matrix}
                          \phi&\phi \\
                          {} & \epsilon
                        \end{matrix}\Big| \lambda^2\right)=1+q+q\phi(-1)\cdot{_2F_1}\left(\begin{matrix}
                          \phi&\phi \\
                          {} & \epsilon
                        \end{matrix}\Big|d^2\right).
  \end{align*}
  Therefore,
\[{_2F_1}(\lambda^2)=\frac{q+1}{q^2}+\frac{q-1}{q}\phi(-1)\;{_2F_1}\left(\left(\frac{1-\lambda}{1+\lambda}\right)^2\right).\]
The other identities follow by observing that Theorem 4.4 of \cite{Greene} indicates that
$_2F_1(d^2)=\phi(-1)\;_2F_1(1-d^2)$ and $\displaystyle _2F_1(d^2)=\phi(1-d^2)\;_2F_1\left(\frac{d^2}{d^2-1}\right)$.
\end{Proof}
The above theorem enables us to obtain new values for $_2F_1(\lambda)$. The values $_2F_1(\pm1)$, $_2F_1(2)$ and $_2F_1(1/2)$ are already known thanks to the works of Greene and Ono, see \cite{Greene, Ono1}. The following result shows how one may obtain new values for $_2F_1(\lambda)$.

\begin{Theorem}
  Let $p\equiv1$ mod $4$ be a prime. Let $a\in\Fp$ be such that $a^2\equiv-1$ mod $p$. One has
  \[_2F_1\left(\frac{4a}{(1+a)^2}\right)=  \frac{2x(-1)^{\frac{x+y+1}{2}}}{p-1}  - \frac{p+1}{p(p-1)},\textrm{\hskip10pt} x^2+y^2=p,\;x \textrm{ is odd}.\]
  Similarly, if $p\equiv \pm1$ mod $8$, we choose $a$ such that $a^2\equiv1/2$ or $2$ mod $8$. One has
  \begin{align*}_2F_1\left(\frac{4a}{(1+a)^2}\right)=\left\{\begin{array}{ll}
  -\frac{p+1}{p(p-1)}& \textrm{if } p\equiv -1 \textrm{ mod }8; \\
\frac{2x(-1)^{\frac{x+y+1}{2}}}{p-1}  - \frac{p+1}{p(p-1)},\textrm{\hskip10pt} x^2+y^2=p,\;x \textrm{ is odd}& \textrm{if } p\equiv 1\textrm{ mod }8. \end{array}\right.\end{align*}
\end{Theorem}
\begin{Proof}
For the first part, the existence of $a$ is guaranteed by the fact that $p\equiv 1$ mod $4$. Now Theorem \ref{thm:transformation} implies that
\begin{align*}
_2F_1\left(\frac{4a}{(1+a)^2}\right)&= \frac{p}{p-1}\left(_2F_1(-1)-\frac{p+1}{p^2}\right)\\
&= \frac{p}{p-1}  \frac{2x(-1)^{\frac{x+y+1}{2}}}{p}  - \frac{p+1}{p(p-1)},\textrm{\hskip10pt} x^2+y^2=p,\;x \textrm{ is odd}
\end{align*}
where the last equality follows from Theorem 2 of \cite{Ono1}.
The second part follows in a similar manner.
\end{Proof}
One may use the transformations given in Theorem \ref{thm:transformation} to evaluate $_2F_1(\lambda)$ for other values of $\lambda$.


\begin{thebibliography}{MM}
\frenchspacing
\renewcommand{\baselinestretch}{1}

\bibitem{BarmanKalita}
R. Barman and G. Kalita,
Certain values of Gaussian hypergeometric series and a family of algebraic curves,
{\em International Journal of Number Theory}, 8 (2012), 945--961.

\bibitem{Barmanalgebraiccurves}
R. Barman and G. Kalita,
Hypergeometric functions and a family of algebraic curves,
{\em Ramanujan Journal}, 28 (2012), 175--185.

\bibitem{BarmanKalitaSaikia}
R. Barman, G. Kalita and N. Saikia,
Hyperelliptic curves and values of Gaussian hypergeometric series,
{\em Archiv der Mathematik}, 102 (2014), 345--355.

\bibitem{Greene}
 J. Greene,
 Hypergeometric functions over finite fields,
 {\em Transactions of the American Mathematical Society},
301 (1987), 77--101.


\bibitem{Gaudry}
P. Gaudry and R. Harley,
Counting points on hyperelliptic curves over finite fields,
In Bosma W. (eds) Algorithmic Number Theory. ANTS 2000. Lecture Notes in Computer Science, vol 1838, 313--332, Springer, Berlin, Heidelberg


\bibitem{Guindy}
A. El-Guindy and K. Ono,
Hasse invariants for the Clausen elliptic curves,
{\em Ramanujan Journal}, 31 (2013), 3--13.


\bibitem{Huff}
G. B. Huff,
Diophantine problems in geometry and elliptic ternary forms,
{\em Duke Mathematical Journal}, 15 (1948), 443--453.

 \bibitem{Joye}
M. Joye, M. Tibouchi, and D. Vergnaud,
Huff's model for elliptic curves,
In {G. Hanrot}, {F. Morain} and E. Thom\'{e}, Eds, Algorithmic
Number Theory (ANTS-IX), LNCS 6197, 234--250, Springer, 2010.

\bibitem{Ono1}
K. Ono,
Values of Gaussian hypergeometric series,
{\em Transactions of the American Mathematical Society}, 350 (1998), 1205--1223

\bibitem{Schoof}
{R. Schoof, Counting points on ellitic curves over finite fields,
{\em Journal de Th\'{e}orie des Nombres de Bordeaux}, 7 (1995), 219--254.}

\bibitem{Sadek1}
M. Sadek and N. El-Sissi,
Edwards curves and Gaussian hypergeometric series,
{\em Journal de Th\'{e}orie des Nombres de Bordeaux}, 28 (2016), 115--124.

\bibitem{Sadek2}
M. Sadek,
Character sums, Gaussian hypergeometric series, and a family of hyperelliptic curves,
{\em International Journal of Number Theory}, 12 (2016), 2173--2187.

\bibitem{S.Z}
{S. Schmitt, H. G. Zimmer, {\em Elliptic Curves: A Computational Approach}, Walter de Gruyter, Berlin  2003.}

\bibitem{Sil}
{J. H. Silverman, {\em The Arithmetic of Elliptic Curves}, Springer, New York 1986.}

\bibitem{Wa}
{L. C. Washington, {\em Elliptic Curves: Number Theory and Cryptography}, 2nd ed., CRC Press, Taylor $\&$ Francis Group, Boca Raton, FL, 2008.}

\bibitem{Wu}
H. Wu and R. Feng,
Elliptic curves in Huff's model,
{\em Wuhan University Journal of Natural Sciences}, 17 (2012), 473--480.

\end{thebibliography}
\end{document}